\documentclass[twoside, 12pt]{amsart}
\usepackage{latexsym}
\usepackage{amssymb,amsmath,amsopn}
\usepackage[dvips]{graphicx}   
\usepackage{color,epsfig}      

\makeatletter
\def\@settitle{\begin{center}%
  \baselineskip14\p@\relax
  \bfseries
  \uppercasenonmath\@title
  \@title
  \ifx\@subtitle\@empty\else
     \\[1ex]\uppercasenonmath\@subtitle
     \footnotesize\mdseries\@subtitle
  \fi
  \end{center}%
}
\def\subtitle#1{\gdef\@subtitle{#1}}
\def\@subtitle{}
\makeatother

\newtheorem{statement}{Statement}
\newtheorem{defi}{Definition}

\begin{document}
\title[]{Elementary Constructions of conic sections}
\subtitle{Constructions without using projective geometry}
\author[\'A. G.Horv\'ath]{\'Akos G.Horv\'ath}
\address {\'A. G.Horv\'ath \\ Department of Algebra and Geometry \\ Institute of Mathematics\\
Budapest University of Technology and Economics\\
M\H uegyetem rkpt. 3, H-1111 Budapest\\
Hungary; \\
ELKH-BME Morphodynamics Research Group, Budapest University of Technology,
M\H uegyetem rakpart 1-3, Budapest 1111, Hungary}
\email{g.horvath.akos@ttk.bme.hu}
\date{}

\subjclass[2010]{00A35, 51M15, 51N05,}

\maketitle
\section{Introduction}
In classical geometry, there is no such well-known and much-studied topic as the construction of conic sections (or briefly conics) from its five points. Its importance in many applications of mechanical engineering, civil engineering and architectural engineering, as well as other applied sciences is clear. The beauty of the topic is that it raises difficult questions that can be approached with basic tools. In this article, we provide constructions (and corresponding theories) that can be taught to high school and university students without knowledge of projective geometry. For this, we recall some important facts about conic sections that can be found in the rich literature. We use the concepts of power of a point on a circle, similarity, orthogonal affinity and inversion. We also mention famous constructions related to our questions. We begin our article at this point, where the standard teaching ends the discussion of conic sections. We therefore assume that the reader knows the basic definitions and constructions of conics, the concepts of focus, axis, tangent, leading circle and leading line.

\section{Definitions}

First, we recall some definitions of conic sections used in our constructions.

\begin{defi}[based on foci]\label{def:fociandsemiaxes}
The \emph{ellipse} is the locus of those points $P$ of the Euclidean plane $E$ for which the sum of distances from two given points $F_1$ and $F_2$ is a constant. The \emph{hyperbola} is the locus of those points $P$ of the plane $E$ for which the absolute value of the difference of its distances from two given points $F_1$ and $F_2$ is a constant. The \emph{parabola} is the locus of those points $P$ of the plane $E$ which are equal distances from a given points $F$ and a given line $l$.
\end{defi}

\begin{defi}[based on a leading circle and one focus]\label{def:focusandleadingcircle}
Let $L$ be a circle in $E$, and ${F}\in E$ be an arbitrary point. The locus of points ${P}\in E$ which are the centre of a circle touches $L$ and go through ${F }$ is called a conic which is defined by its leading circle $L$ and its focus $F$. It is an \emph{ellipse} if $F$ is in the interior of $L$, is a \emph{hyperbola} if $F$ lies in the exterior of $L$, and it is a \emph{parabola} if the radius of $L$ is infinite, so $L$ is a line.
\end{defi}

\begin{defi}[based on a leading line and a focus]\label{def:focusandleadingline}
A conic is the locus of points in $E$, for which the ratio of their distances from a given point $F$ and a given line $l$ is a constant. It is an ellipse, parabola or hyperbola if this ratio is less, is equal or is greater than 1, respectively.
\end{defi}

The equivalence of these definitions is well-known and can be found in all text-book on conics. For the non-familiar reader, I propose the book \cite{glaeser-stachel-odehnal}. The following metric definition seems to be less well-known. Budden observed in \cite{budden} that the affine coordinates with respect to such an affine coordinate system which axes are a pair of conjugate diameters give also well-usable metric definitions on conics. He proved (with elementary methods) the following statement:

\begin{statement}\label{def:budden}
The locus of points $P$ which satisfies the equation
$$
\frac{PV^2}{VQ\cdot VQ'}=\nu
$$
a constant in magnitude and sign, $QQ'$ being fixed real points, and $PV$ parallel to a fixed direction meeting $QQ'$ in $V$, is a hyperbola or an ellipse, according as $\nu$ is positive or negative.

The locus of points
$$
PV^2=\nu\cdot QV,
$$
where $\nu$ is a constant positive length, and $PV$, parallel to a fixed direction, meets fixed line through a fixed point $Q$ in $V$, is one or other of two parabolas, according to the direction of $QV$ which is considered positive.
\end{statement}

\begin{figure}[h]
\includegraphics[scale=1]{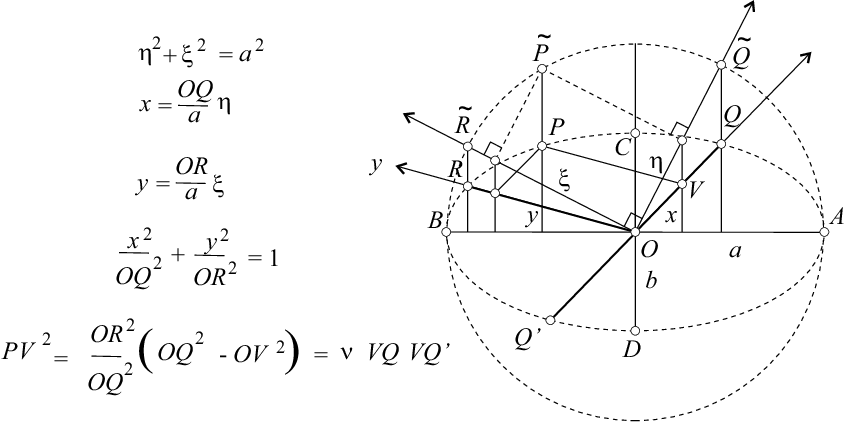}
\caption{Affine equation of an ellipse and the Budden's definition. }
\label{fig:affineequation}
\end{figure}

Clearly, this statements can be considered alternative definitions for conic sections. In the case of ellipse, we can check it in Figure \ref{fig:affineequation}. In Figure \ref{fig:affineequation}., we denoted by $(\eta, \xi)$ the coordinates of the point $\tilde{P}$, with respect to the cartesian system $\{O\tilde{Q}, O\tilde{R}\}$. Moreover, let $(x,y)$ be the affine coordinates of the point $P$ with respect to the conjugate diameters $\{OQ,OR\}$. Then we can express $PV=y$ by the signed distances of $VQ$, and $VQ'$ as follows
$$
PV^2=\frac{OR^2}{OQ^2}\left(OQ^2-OV^2\right)=\nu(OQ-OV)(OQ+OV)=\nu QV\cdot Q'V,
$$
where $\nu=\frac{OR^2}{OQ^2}$.

\section{Inversion in action}

Consider the point $O$ of the plane and project the point $P$ of the plane onto the point $P'$ of the half-line $\overrightarrow{OP}$ to which $OP\cdot OP' =r^2$ is satisfied by the given and fixed non-zero real number $r$. Then we say that $P'$ is the image of $P$ under the inversion is defined by the circle with center $O$ and radius $r$. Inversion maps circles or lines onto other circles or lines, and this preserves the angles between the crossing curves. With the help of inversion, problems about circles can be transformed into problems about other lines. This subsection contains problems on conic sections can be solved with inversion.

\subsection{Conic section with given focus and three points}

\begin{figure}[ht]
\centerline{\includegraphics[scale=0.8]{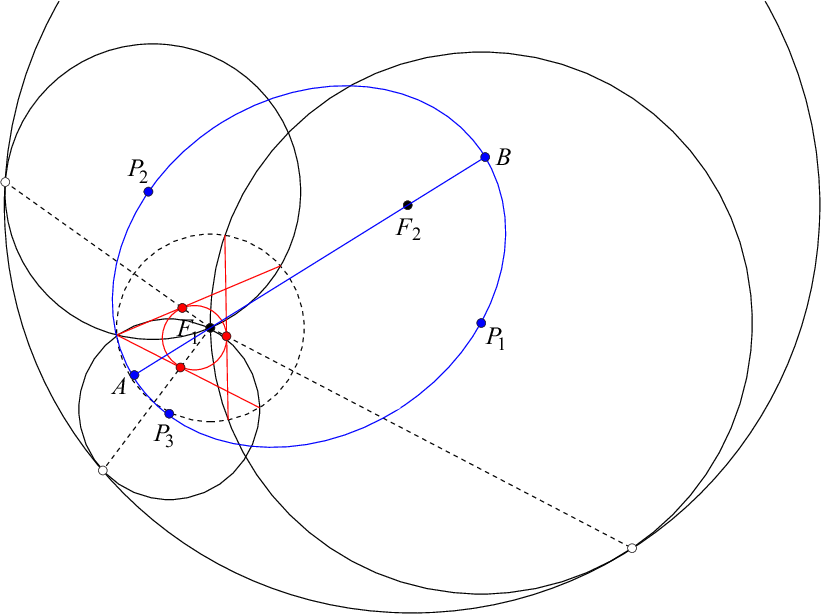}}
\caption{The construction of the second focus if given a focus and three points of the conic section.}
\label{fig:threepoints}
\end{figure}

\emph{One focus and three points of the conic section are given. Let's construct the essential data of the conic section!}
The solution to the problem is based on Definition \ref{def:focusandleadingcircle} and the inversion properties. Rewriting the example with this definition reads: \emph{Given three circles through a point. Let's draw a fourth circle that touches all three.} The inversion of the common point as the center turns the circles into pairs of intersecting lines. The four circles touching all three lines are the inverse image of the four possible solutions to the problem. The sought circles are determined by the inverse images of the points of contact on the inverse lines.

\begin{remark}Defining a circle that touches three given circles is usually a more complicated task. In this case, where each circle lies on the outside of the other two, Gergonne gives a possible solution (see e.g. \cite{dorrie} or \cite{gho_malfatti}):
\emph{Draw the point $P$ of power of the given circles and an axis of similitude of certain three centres of similitude. The poles of the axis of similitude with respect to the given circles are points of the rays from $P$ passing through the touching points of a sought circle and one of the three given ones.}

\begin{figure}[htbp]
\centerline{\includegraphics[scale=1]{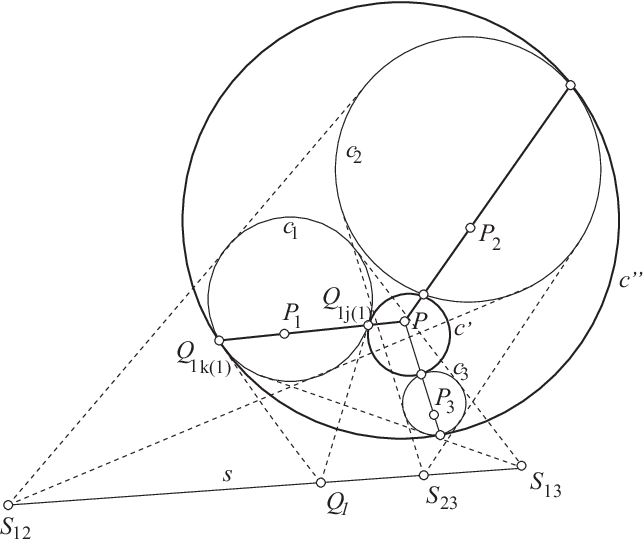}}
\caption{The construction of Gergonne}
\label{fig:gergonne}
\end{figure}.

The proof of this construction follows from the pole-polar theorem and from the fact that the inversion for circles perpendicular to the given circles with center $P$ leaves the given circles invariant. It follows that you exchange the solution that touches the circles from the outside with the one that contains the three given circles. Therefore, as we can see in Figure \ref{fig:gergonne} the points $P,Q_{1,j_1},Q_{1,k_1}$ are collinear, and the pole of the line $Q_l$ matches $s$. (We note that the similarity center $S_{1,2}$ is the center of the inversion, which swaps the circles $c_1$ and $c_2$ and leaves the tangent circles $c'$ and $c''$ invariant. The $c' $ and the power of the circles $c''$ with respect to the point $S_{1,2}$ are equal. Similar properties also apply to $S_{1,3}$, proving this claim.)
\end{remark}

\subsection{Conic section with a focus, two points and a tangent }

\emph{One focus, two points and a tangent of the conic section are given. Let's construct the essential data of the conic section!}. Paraphrasing the task by Definition \ref{def:focusandleadingcircle}, we have to construct that circle which goes through a point (which is the reflected image of the focus to the touching line) and touching two circles (with center of the given points through the focus). Clearly, an inversion on the focus transforms the task to construct circle through a point and touching two lines intersecting in the second common point of the given circles. This latter problem can be solved by similarity easily. Applying the inversion again, we get that leading circle of the conic which centre is the other focus. If the inverted image is a line (its centre is a point at infinity) we have got the leading line of a parabola.

\begin{figure}[ht]
\centerline{\includegraphics[scale=0.8]{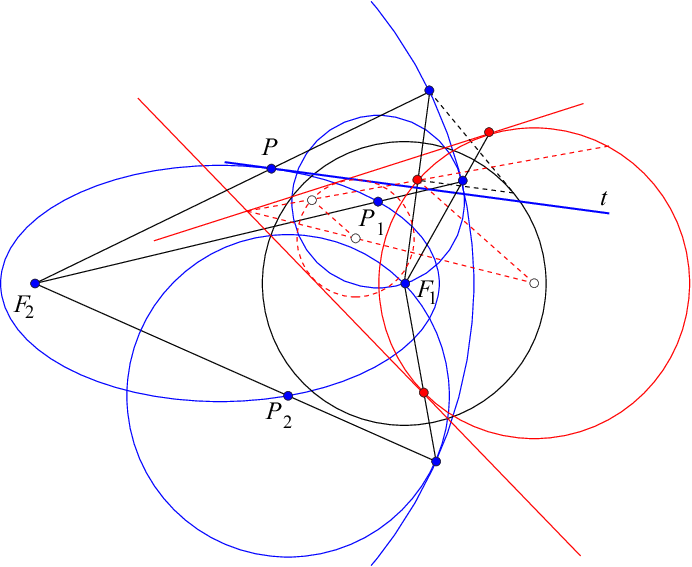}}
\caption{The construction of the second focus if given two points and a tangent of the conic section. Inverted elements are marked in red.}
\label{fig:twopointsandatangent}
\end{figure}

\begin{remark} This construction leads to a group of geometric construction in which we have to construct a circle which touches some other circles or lines. One of these problems is: \emph{Given a point, a line and a circle. Construct a circle touching the given line and circle through the given point.}
Observe that an inversion with respect to a circle which centre is on the line and goes through the point, transforms the problem to the following: Given two circle and a point which lies on the first one. Construct a circle which in the given point touches the first circle and also touches the second one. That is we should construct a circle which centre is in a line $t$ goes through a point $P$ of this line and touches a given circle $k$ with centre $K$. If the radius of $k$ is $r$ the distance of $P$ from the centre $C$ of the solution is $d$ then $|CK|=d\pm r$. Therefore, for $Q\in t$ $|QP|=r$, $C$ falls at the same distance from $K$ and $Q$. The bisector of the $KQ$ segment intersects $t$ in the sought point $C$.
\end{remark}

Note that by dilation and inversion any such problem can be simplified to one of the following:
\begin{itemize}
\item \emph{Given three lines. Construct a circle touching the given lines.}
\item \emph{Given a point, and two lines. Construct a circle through the given point touching the given lines.}
\item \emph{Given a point, a line and a circle. Construct a circle through the given point touching the given line and circle.}
\item \emph{Given a point, and two circles. Construct a circle through the given point touching the given circles.}
\end{itemize}
The solution of the first one elementary. The easy solution of the second problem uses similarity. As we saw the solution of the third problem can be transformed by a special case of the fourth problem. We remark that a construction of the last problem solves the problem of the preceding paragraph substituting the Gergonne's construction. Assume that the given circles are denoted by $(K_i,r_i)$, where their centres are $K_i$ and radiuses are $r_i$. If $r_1\leq r_2\leq r_3$ holds then the centre $K$ of the sought circle of radius $r$ is also the centre of a circle of radius $r\pm r_1$ which goes through $K_1$ and touches the circles $(K_2,r_2\pm r_1)$ and $(K_3,r_3\pm r_1)$, respectively. On the other hand, this problem applying an inversion in a point of $(K_2,r_2\pm r_1)$ transforms to the third problem which we solved in this subsection. Therefore, with inversion, we can solve the problem of construction circles touching three given circles without using projective geometry tools.

\section{The role of the affinities}

In the following  constructions, affine mappings play an important role.

\subsection{Ellipse from the line of the axes and two known points.}\label{subsec:ellipse}

\begin{figure}[ht]
\includegraphics[scale=0.7]{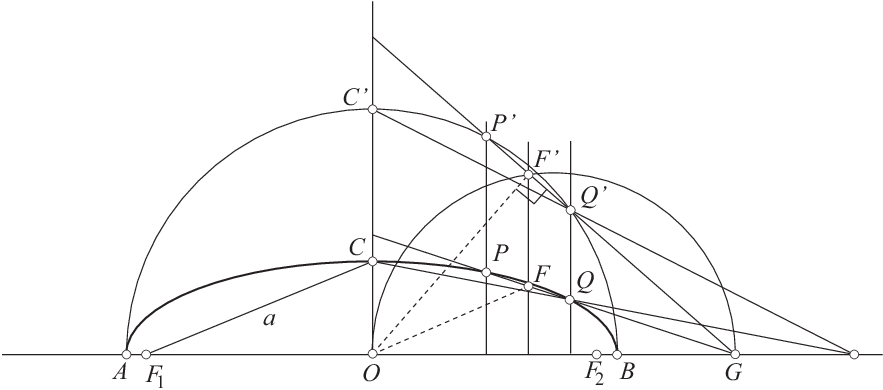}
\caption{Construction of an ellipse from its axes and two points.}\label{fig:ellipsefromaxes}
\end{figure}

Let $P$ and $Q$ be the given points and denote by $t$ one of the given axe through the centre $O$ of the ellipse. Let $F$ be the midpoint of the segment $PQ$ and denote by $G$ the intersection point of the line $t$ with the line $PQ$. The circle with diameter $OG$ meets the vertical line through $F$ in the points $F'$. The affinity which axis is $t$ and sends the points $F$ to $F'$ and $P$ to $P'$ also sends $Q$ to $Q'$. Now the length of the segment $OP'$ is equal to the length of the segment $OQ'$ implying that the image of the searched ellipse at this affinity is a circle. The radius of this circle is the length of the major semi-axis and we can get easily also the endpoint of the minor semi-axis using this connections. (See Figure \ref{fig:ellipsefromaxes}.) The foci are obtained from the axes with the well-known manner.

\subsection{Hyperbola from two asymptotes and a known point}\label{subsec:hyperbola}

\begin{figure}[ht]
\includegraphics[scale=0.7]{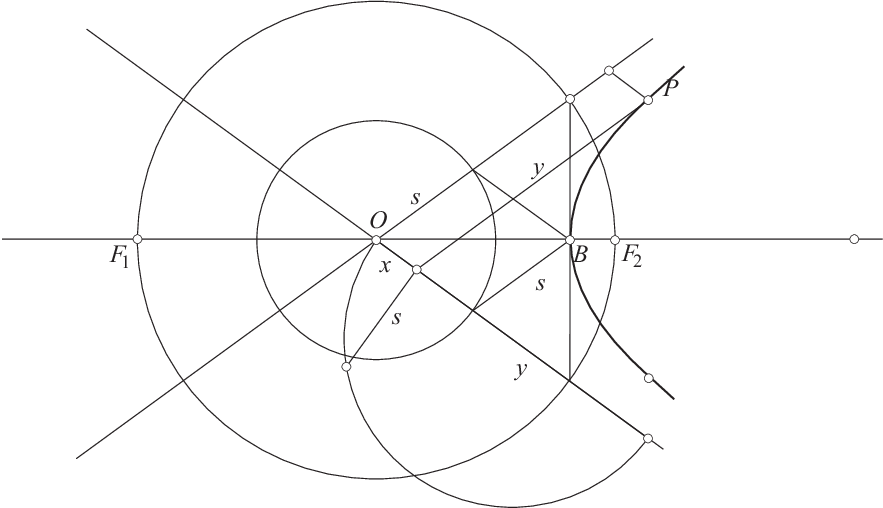}
\caption{Construction of hyperbola from its asymptotes and one of its points.}\label{fig:hyperbolafromasymptotes}
\end{figure}

First we construct those vertex $B$ of the real axis which bisects the angle domain of the asymptotes containing the given point $P$. Since the product of the affine distances $x$ and $y$ of the point $P$ from the asymptotes is independent of the position of $P$ and for the point $B$ this distances are equal to each other, we get the common value is the geometric mean $s$ of $x$ and $y$. The tangent $b$ at the vertex $B$ intersects the asymptotes two points which are lying on the circle with center $O$ and through the foci. (See in Figure \ref{fig:hyperbolafromasymptotes}.)

\subsection{Parabola from its axis and two of its points}\label{subsec:parabola}

\begin{figure}[ht]
\includegraphics[scale=1]{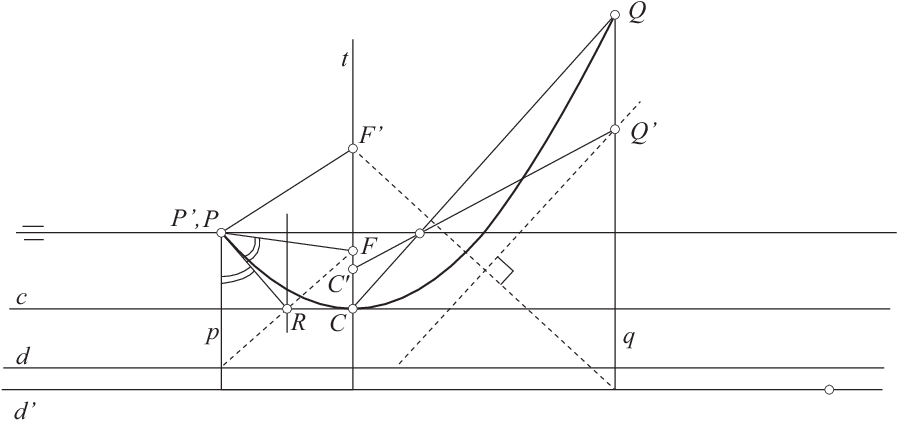}
\caption{Construction of parabola from its axis and two of its points.}\label{fig:parabolafromaxis}
\end{figure}

We use the fact that the affine copy of a parabola is also a parabola. Clearly, if the axis of an orthogonal affinity is parallel to the direktrix then the axis of the image is the same as the axis of the original one. This means that the axis of the symmetry is unchanged and the image of the vertex of the parabola is the vertex of the image parabola. On the axis $t$ choose a point $F'$ and determine the direktrix $d'$ such that one of the given points say $P$ will be a point of the parabola $(F',d')$. Denote by $C'$ its vertex. On that diameter $q$ of this parabola which contains the other given point $Q$ determine the unique point of $(F',d')$. (Of course, if intersects the line $q$ the line $d'$ at the point $E$, the orthogonal bisector of the segment $F'E$ intersects the point $Q'$ from the line $q$.) Consider that orthogonal affinity which axis is the line $s$ through the point $P$ and perpendicular to the axis $t$ of $(F',d')$ and sends the point $Q'$ to the point $Q$. This affinity sends $P$ to $P$ and $C'$ to a point $C$ of $t$ which is the vertex of the image parabola. Let $p$ denote the common diameter of the two parabolas passing through point $P$, and let $c$ denote the tangent of the image parabola at the vertex $C$. Finally, let $R$ be the intersection point of the tangent $c$ with the middle parallel of the parallel lines $p$ and $t$. It is known that $RP$ is the tangent to the image parabola at point $P$, and the reflection of the line $p$ to the tangent $RP$ intersects the axis $t$ at the focus $F$. Now we also get the directrix $d$ from the points $F$, $C$. (See in Figure \ref{fig:parabolafromaxis}.)

\subsection{Central conic from a diameter, a point and the direction of the conjugate diameters}

\begin{figure}[ht]
\includegraphics[scale=1]{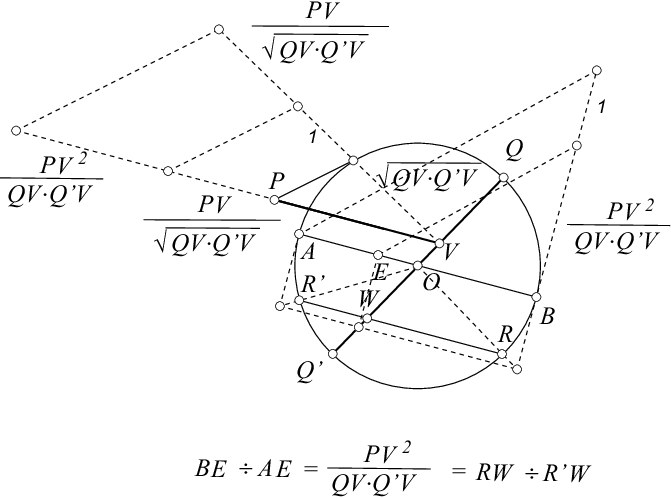}
\caption{Construction of the ratio $\frac{PV^2}{VQ\cdot VQ'}$ and the chord $RR'$ in the Thales circle. }
\label{fig:buddenratioconstr}
\end{figure}

Using Definition \ref{def:budden} we have the ratio
$$
\nu=\frac{PV^2}{VQ\cdot VQ'},
$$
which can be constructed easily (see Figure \ref{fig:buddenratioconstr} ). If $QQ'$ cuts a chord of the Thales circle of $QQ'$ parallel to $PV$ in ratio $R'W:RW=\nu$ then $R$ is a point of the conics and the bisecting line of $QR$ is an axis of the conic. In fact, since $QRQ'R'$ are concyclic $WQ\cdot WQ'=WR\cdot WR'$, hence
$$
\frac{RW^2}{WQ\cdot WQ'}=\frac{WR^2}{WR\cdot WR'}=\nu,
$$
implying that $R$ is on the conic defined by $P,V,Q,Q'$ and $\nu$. On the other hand $QR$ perpendicular to $RQ'$ by Thales theorem hence $QR$ and $RQ'$ are parallel to the two axes, respectively. (Two conics with a common center have at most four points in common. Therefore, a circle, which has its center in common with a central conic, intersects it at the vertices of a rectangle. It is possible that the four points degenerate into two, in this case they are on one axis of the cone.) $QR$ and the parallel to $PV$ from $Q$ intersect the bisector line of $QR$ in $M$ and $T$, respectively. Let $A$ be that point of this axis for which holds $OA^2=OM\cdot OT$. Then $A$ is in the conic and we have found the vertices of the conic in this axis. In fact, we have
$$
\frac{OM}{OA}=\frac{OA}{OT}=\frac{EO}{OE},
$$
hence $EM$ is parallel to $QA$. So $EM\cdot OQ=QA\cdot OE$ from which by adding $EM\cdot OE$ we get that
$$
\frac{EK}{EA}=\frac{EM}{EM+QA}=\frac{OE}{OE+OQ}=\frac{OE}{Q'E}=\frac{EA}{EL}.
$$
So $EK\cdot EL=EA^2=AE^2$ therefore
$$
\frac{AE^2}{EQ\cdot EQ'}=\frac{EK\cdot EL}{EQ\cdot EQ'}=\frac{WR\cdot WR}{WQ\cdot WQ'}=\nu,
$$
as we stated.

\begin{figure}[ht]
\includegraphics[scale=1]{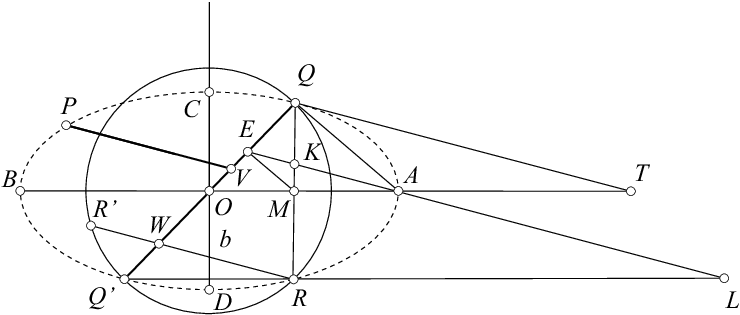}
\caption{Budden's construction.}\label{fig:buddenconstr}
\end{figure}

\section{Pascal line and conics from five points.}

\subsection{Pascal's theorem on hexagons inscribed into a conic}

\begin{figure}[ht]
\includegraphics[scale=0.8]{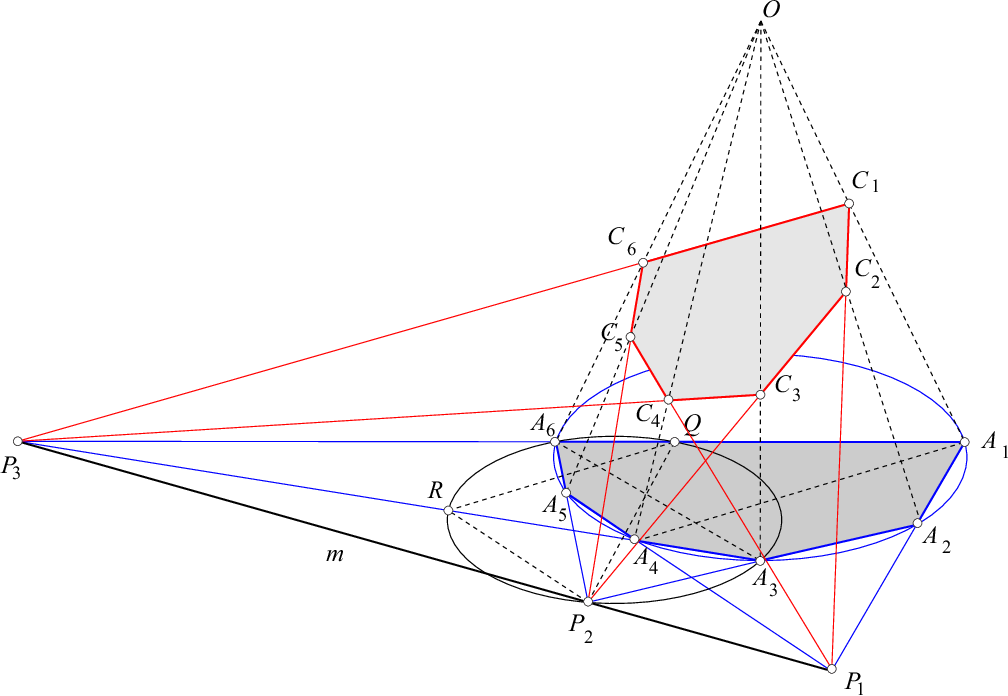}
\caption{Elementary proof of Pascal's theorem for conic section.}\label{fig:centralprojectofcircle}
\end{figure}

In this section we construct the important data of a conic from its five points. We use the affine constructions of the previous paragraph and a new tools the Pascal line of the conics. Pascal theorem says that the opposite sides of a hexagon inscribed in a conic meets on a line the so-called Pascal line. Clearly, the position of the line depends on the order of the points, so for six points we can associate sixty lines in general. On the other hand, in our construction doesn't important which Pascal line helps us. (See more on this problem in \cite{glaeser-stachel-odehnal}.) The most simple Euclidean proof of Pascal's circle theorem can be found in paper \cite{yzeren}, it is based on similarity. We reproduce this nice proof and extend it to the general case. (See the notation in Figure \ref{fig:centralprojectofcircle}.) Consider the hexagon $A_1A_2\ldots A_6$ inscribed in the blue circle. The corresponding opposite sides intersect in pairs at $P_1$, $P_2$ and $P_3$. Construct the dashed circle containing the points $A_3,A_6$ and $P_2$. This intersects the lines $P_3A_6$ and $P_3A_4$ in the points $R$ and $Q$, respectively. Using equal angles, we can prove easily that the segments $RQ$, $RP_2$ and $P_2Q$ are parallel to the segments $A_1A_4$, $A_4P_1$ and $P_1A_1$, respectively. This implies that the central similarity with center $P_3$ and sending $Q$ to $A_1$, also sends $R$ and $P_2$ to $A_4$ and $P_1$, respectively. Hence $P_1$, $P_2$ and $P_3$ are collinear. Consider a circular cone above the $\alpha$ plane of the above circles, which passes through the circumscribed circle of the hexagon $A_1\ldots A_6$. Take a plane $\beta$ which intersects $\alpha$ at the common line $m$ of the points $P_i$. $\beta$ intersects the generator $OA_i$ of the cone in the points $C_i$, determining a new hexagon $C_1C_2\ldots C_6$ inscribed in the corresponding conic sections of the cone. Clearly, this hexagon has the Pascal's property with the Pascal line $m$. For every conic section, we can get a similar one, with an appropriately chosen $\beta$ plane, and every hexagon of type $C_1C_2\ldots C_6$ inscribed in this conic slice is an image of a hexagon $A_1A_2\ldots A_6$ inscribed in the blue circle. So, we get that every hexagon inscribed in a conic section has the Pascal property.

\subsection{Construction of a conic section from five points}

\begin{figure}[ht]
\includegraphics[scale=0.8]{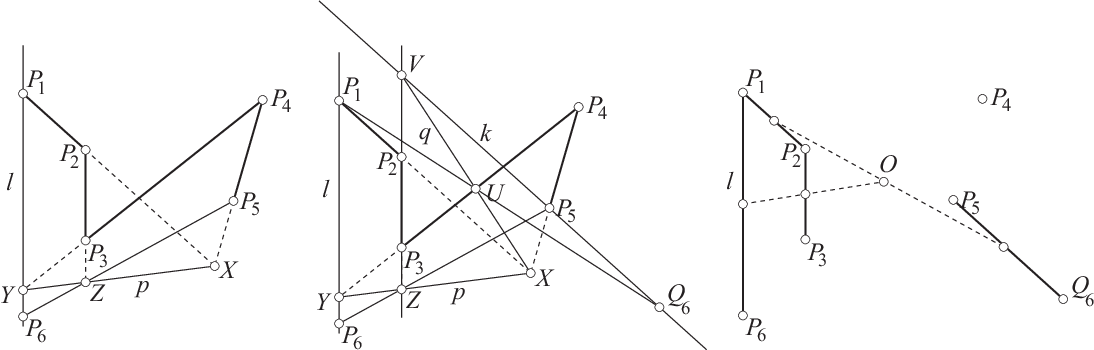}
\caption{Construction of parallel chords and the center of the conic.}\label{fig:paralelchords}
\end{figure}

Let $P_1$, $P_2,\ldots P_5$ be five points on the plane. We construct the essential data of that conic section defined by these points. First we determine two pairs of parallel chords of the conic. To this consider a line $l$ through the point $P_1$ and parallel to the segment $P_2P_3$ and determine the second point of intersection of $l$ with the conic. The Pascal line $p=(XYZ)$ is determined by the points $X=P_1P_2\cap P_4P_5$ and $Y=P_3P_4\cap l$. The line $P_2P_3$ meets $XY=p$ at the point $Z$ and $ZP_5$ intersects $l$ in $P_6$. (See the left figure in Figure \ref{fig:paralelchords}. ) Now $P_1P_6$ and $P_2P_3$ are parallel chords of the conic. On a similar way we can construct the point $Q_6$ with the property that $P_1P_2$ parallel to $P_5Q_6$. (See the middle picture in Figure \ref{fig:paralelchords}.)

Since the midpoints of the parallel chords lie on one diameter of the conic, the line connecting the midpoints of the parallel chords contains the center of the conic. (See the right picture in Figure \ref{fig:paralelchords}.) From a Euclidean point of view, we have two different options.

\begin{figure}[ht]
\includegraphics[scale=0.8]{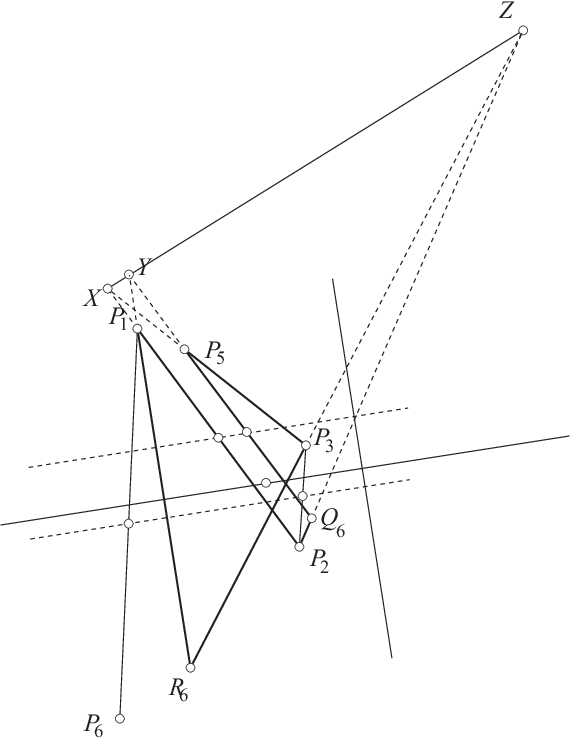}
\caption{Construction from the direction of the diameters the axis of the parabola.}\label{fig:axisoftheparabola}
\end{figure}

\begin{itemize}
\item In the first case, the diameters are parallel, which means that the center lies at infinity, i.e. the conic section is a parabola. In this case, we know the direction of the axis and can determine the point $R_6$, which is the second point of the parabola on the line that is perpendicular to the diameters and passes through one of the points (say, the point marked $P_1$). The perpendicular bisector of the section $P_1R_6$ is the axis of the parabola, and the focus and direction can be constructed with the result of Subsection (\ref{subsec:parabola}). (See the construction of the axis in Figure \ref{fig:axisoftheparabola}.)

\item In the second case, we must distinguish between the case of a hyperbola and an ellipse. Two pairs of conjugate diameters are available at the $O$ center, so if we take a line that intersects these lines in four points, we have two possibilities.
    \begin{itemize}
    \item In the first, the corresponding pairs of points cross each other, creating a negative cross-product. These is the case of the ellipse. This situation can be solved with the method of Subsection (\ref{subsec:ellipse}).
\begin{figure}[ht]
\includegraphics[scale=0.8]{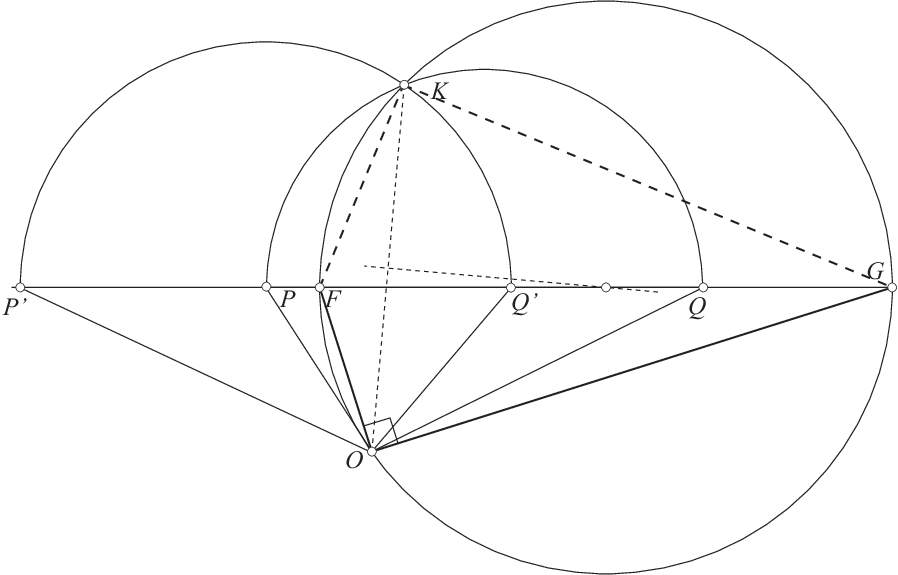}
\caption{Construction from the direction of the diameters of the ellipse its axes.}\label{fig:axesoftheellipse}
\end{figure}
        First, determine the axes of the ellipse as follows (see Figure \ref{fig:axesoftheellipse}): Consider a skew-affinity that sends the two pairs of conjugate diameters into pairs with orthogonal elements. A circle passing through the respective centers $O$ and $K$, centered on the axis of affinity, intersects this line at two points ($F$ and $G$ in the figure), which are the points of the axes of the ellipse.

        In the second step apply the construction of Subsection (\ref{subsec:ellipse}) to get the required data of the ellipse.

    \item If the two pairs of points do not intersect, then we get the case of a hyperbola. First we define the asymptotes and then apply the construction of subsection \ref{subsec:hyperbola}. Note that we can determine the asymptotes of the hyperbola from the points $P$, $Q$ and $R$ and the center $O$ (see Figure \ref{fig:asymptotesofthehyperbola}).
\begin{figure}[ht]
\includegraphics[scale=0.8]{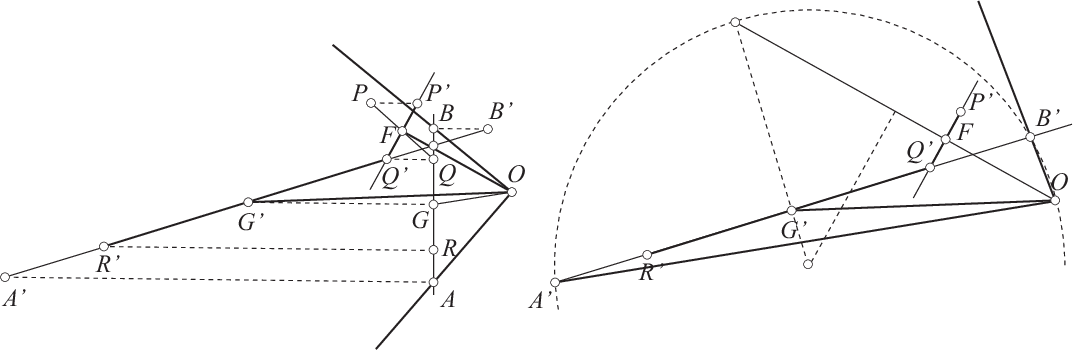}
\caption{Construction from the direction of the diameters of the hyperbola its asymptotes.}\label{fig:asymptotesofthehyperbola}
\end{figure}
         In fact, by applying skew-affinity, we can map the segment $PQ$ to the segment $P'Q'$, which is perpendicular to the line $OF$, where $F$ is the midpoint of $PQ$ (see the left image in Figure \ref{fig:asymptotesofthehyperbola}). This affinity maps $R$ to $R'$ and the midpoint $G$ of the segment $QR$ to the midpoint $G'$ of the segment $Q'R'$ . The asymptotes of the new hyperbola defined by the points $P', Q', R'$ and the center $O$ can be obtained using the method of the right-hand image of Figure \ref{fig:asymptotesofthehyperbola}. In this picture, we use the fact that the midpoint $G'$ of the segment $Q'R'$ is also the midpoint of the segment $A'B'$, where the sought asymptotes are $OA'$ and $OB'$. Therefore, we know the vertex $O$ of the triangle $OA'B'$, the angle bisector at this vertex, which is $OF$, and the median $OG'$ of the triangle. It is known that the segment bisector perpendicular to the side $A'B'$ intersects the angle bisector $OF$ at a point of the circumscribed circle of the triangle, so we can construct the vertices $A'$ and $B'$ as the intersection of the circumscribed circle with the line $A'B'$. Finally, using the previous affinity again, we get the points $A$ and $B$ of the asymptotes from the points $A'$ and $B'$ (see the left image of Figure \ref{fig:asymptotesofthehyperbola}.
    \end{itemize}

\end{itemize}

\end{document}